\documentclass[12pt]{amsart}

\usepackage{graphicx,amssymb,latexsym,amsfonts,txfonts,amsmath,amsthm}
\usepackage{pdfsync,color,tabularx,rotating}
\usepackage[all,cmtip]{xy}
\usepackage{url}
\usepackage{tikz}
\usepackage{amssymb}
\usepackage{comment}

\advance\hoffset by -0.9 truecm

\newtheorem{theo}{Theorem}[section]

\newtheorem{conj}[theo]{Conjecture}

\theoremstyle{remark}
\newtheorem{rema}[theo]{\bf Remark}

\theoremstyle{remark}
\newtheorem{exam}[theo]{\bf Example}

\theoremstyle{remark}

\theoremstyle{remark}

\begin{document}

\title{Orders of simple groups and the Bateman--Horn Conjecture}

\author{Gareth A.~Jones and Alexander K.~Zvonkin}

\address{School of Mathematical Sciences, University of Southampton, Southampton SO17 1BJ, UK}
\email{G.A.Jones@maths.soton.ac.uk}

\address{LaBRI, Universit\'e de Bordeaux, 351 Cours de la Lib\'eration, F-33405 Talence Cedex, France}

\email{zvonkin@labri.fr}


\begin{abstract}
We use the Bateman--Horn Conjecture from number theory to give strong evidence of a positive answer to Peter Neumann's question, whether there are infinitely many  simple groups of order a product of six primes. (Those with fewer than six were classified by Burnside, Frobenius and H\"older in the 1890s.) The groups satisfying this condition are ${\rm PSL}_2(8)$,  ${\rm PSL}_2(9)$ and ${\rm PSL}_2(p)$ for primes $p$ such that $p^2-1$ has just six prime factors. The conjecture suggests that there are infinitely many such primes, by providing heuristic estimates for their distribution which agree closely with evidence from computer searches. We also briefly discuss the applications of this conjecture to other problems in group theory, such as the classifications of permutation groups and of linear groups of prime degree, the structure of the power graph of a finite simple group, and the construction of highly symmetric block designs.
\end{abstract}

\maketitle


\section{Introduction}\label{sec:intro}

In the 1890s, Burnside~\cite{Bur94, Bur95}, Frobenius~\cite{Fro93} and H\"older~\cite{Hol} classified those non-abelian finite simple groups $G$ which have a number-theoretically small order, in the sense that $|G|$ is a product of relatively few primes, counting repetitions. Between them, they showed that there are only four such groups with at most five primes, namely the groups ${\rm PSL}_2(p)$ for $p=5$, $7$, $11$ and $13$, of orders
\[2^2\cdot 3\cdot 5,\quad 2^3\cdot 3\cdot 7,\quad 2^2\cdot 3\cdot 5\cdot 11\quad\hbox{and}\quad 2^2\cdot 3\cdot 5\cdot 13.\]
(Recall that ${\rm PSL}_2(q)$ (or $L_2(q)$ in ATLAS~\cite{CCNPW} notation) has order $q(q^2-1)/d$ where $d=\gcd(q-1,2)$, and that ${\rm PSL}_2(4)\cong{\rm PSL}_2(5)\cong{\rm A}_5$.)

Since the non-abelian finite simple groups are now classified, it is natural to ask which of them have order a product of six primes. It is straightforward to inspect the orders of the simple groups in such sources as~\cite{CCNPW} or~\cite{Wil}, and to see that the only possibilities are the groups ${\rm PSL}_2(8)$ and ${\rm PSL}_2(9)$ ($\cong{\rm A}_6$) of orders
\[2^3\cdot 3^2\cdot 7\quad\hbox{and}\quad 2^3\cdot 3^2\cdot 5,\]
and ${\rm PSL}_2(p)$, of order $p(p^2-1)/2$, for primes $p>13$ such that $p^2-1$ is a product of six primes. 
For instance, see Example~\ref{ex:unitary} in Section~\ref{sec:more} for the elimination 
of small unitary groups.

Solomon, in his historical survey~\cite{Sol} of the classification of finite simple groups, notes that Peter Neumann, while editing Burnside's collected works, asked whether there are infinitely many primes $p$ with this property. We do not have a definite answer to this question, but we have used the Bateman--Horn Conjecture to give what we feel is strong evidence that the answer is positive.


\section{Factorisations}\label{sec:factors}

For a natural number $n\in{\mathbb N}$ with prime power factorisation
\[n=\prod_{i=1}^kp_i^{e_i}\quad (p_i\;\hbox{prime},\; e_i\ge 1),\]
let
\[\Omega(n):=\sum_{i=1}^ke_i\]
denote the total number of prime factors of $n$, counting repetitions, and for a finite group $G$ let
\[\Omega(G):=\Omega(|G|).\]
For each $m\in{\mathbb N}$ let ${\mathcal S}_m$ be the set of all non-abelian simple groups $G$ such that $\Omega(G)=m$. As explained in the Introduction, ${\mathcal S}_6$ consists of the groups ${\rm PSL}_2(8)$ and ${\rm PSL}_2(9)$, together with the groups ${\rm PSL}_2(p)$ for those primes $p>13$ such that $\Omega(p^2-1)=6$.

For any integer $p$ we have
\[\Omega(p^2-1)=\Omega(p-1)+\Omega(p+1).\]
Now assume that $p$ is prime and $p>13$. Then $p\pm 1$ are both even, one of them is divisible 
by~$4$, and one of them is divisible by $3$. This gives us four prime factors contributing to $\Omega(p^2-1)$, and each of $p\pm 1$ must contribute at least one more prime factor since otherwise $p\le 13$. Thus $\Omega(p^2-1)\ge 6$, and this lower bound is attained if and only if $p$ satisfies one of the following conditions, depending on which of $p\pm 1$ is divisible by $3$ or by $4$:
\begin{itemize}
\item[a)] $p-1=4r$ and $p+1=6s$,
\item[b)] $p-1=6r$ and $p+1=4s$,
\item[c)] $p-1=2r$ and $p+1=12s$,
\item[d)] $p-1=12r$ and $p+1=2s$,
\end{itemize}
where $r$ and $s$ are primes. All four of these cases are possible, for example with $p=29$, $19$, $23$ and $37$ respectively, but the question is whether any of them yields infinitely many primes $p$.

In case~(a) we have $p\equiv 1$ mod~$(4)$ and $p\equiv -1$ mod~$(6)$, which are equivalent to $p\equiv 5$ mod~$(12)$. If we put $p=12t+5$ where $t\in{\mathbb N}$ then $r=3t+1$ and $s=2t+1$. The question is now whether these three polynomials
\[f_i(t)=12t+5,\quad 3t+1\quad\hbox{and}\quad 2t+1\quad (i=1, 2, 3)\]
can simultaneously take prime values for infinitely many $t\in{\mathbb N}$. The situation is similar in the other three cases, the only difference being that the triples of polynomials $f_i(t)$ are $12t+7$, $2t+1$ and $3t+2$ in case~(b), $12t-1$, $6t-1$ and $t$ in case~(c), and  $12t+1$, $t$ and $6t+1$ in case~(d).


\section{Number-theoretic conjectures}\label{sec:conjs}

These four cases are instances of a much more general problem in number theory, namely whether a given set of polynomials $f_1(t),\ldots, f_k(t)\in{\mathbb Z}[t]$ can simultaneously take prime values for infinitely many $t\in{\mathbb N}$. Bunyakovsky~\cite{Bun} considered the case $k=1$ in 1857. The following conditions are obviously necessary for a single polynomial $f(t)\in{\mathbb Z}[t]$ to take prime values for infinitely many $t\in{\mathbb N}$:
\begin{itemize}
\item[(1)] it has a positive leading coefficient,
\item[(2)] it is irreducible in ${\mathbb Z}[t]$, and
\item[(3)] it is not identically zero modulo any prime.
\end{itemize}
Bunyakovsky conjectured that these conditions are also sufficient. For example, they are satisfied by the polynomial $t^2+1$, the subject of the Euler--Landau Conjecture on primes of this form. The Bunyakovsky Conjecture has been proved only in the case where $\deg f=1$: this is simply a reformulation of Dirichlet's Theorem on primes in an arithmetic progression.

Our cases (a) to (d) are instances of Dickson's Conjecture, which is that polynomials $f_i(t)$ \linebreak ($i=1,\ldots, k$) of degree $\deg f_i=1$ simultaneously take prime values for infinitely many $t\in{\mathbb N}$ if and only if they all satisfy the first two Bunyakovsky conditions and their product satisfies the third. Particular cases include the twin primes and Sophie German primes conjectures, with $f_i(t)=t$, $t+2$ and $t$, $2t+1$ respectively. Again, this conjecture has been proved only in the case $k=1$, whereas our six primes problem has $k=3$.

In 1957 Schinzel's Hypothesis~\cite{SS} generalised both the Bunyakovsky and Dickson Conjectures by removing the condition $\deg f_i=1$ from the latter. In 1962 Bateman and Horn~\cite{BH}, extending earlier work of Hardy and Littlewood~\cite{HL} on the twin primes and other related conjectures, proposed an asymptotic estimate $E(x)$ for the number $Q(x)$ of $t\in{\mathbb N}$ with $t\le x$ such that $f_i(t)$ is prime for $i=1,\ldots,k$. The Bateman--Horn Conjecture (BHC) asserts that 
\begin{equation}\label{eqn:E}
Q(x)\sim E(x):=C\int_a^x\frac{dt}{\prod_{i=}^k\ln f_i(t)}\quad\hbox{as}\quad x\to+\infty
\end{equation}
where $C$, known as a {\sl Hardy--Littlewood constant}, is given by
\begin{equation}\label{eqn:C}
C=C(f_1,\ldots,f_k):=\prod_{r\;{\rm prime}}\left(1-\frac{1}{r}\right)^{-k}\left(1-\frac{\omega_f(r)}{r}\right),
\end{equation}
with $\omega_f(r)$ denoting the number of roots of $f:=f_1\ldots f_k$ mod~ $(r)$, and $a$
in (\ref{eqn:E}) chosen to avoid singularities of the integral, where some $f_i(t)=1$. In the next section we will give a short heuristic argument to explain these formulae, but as yet there is no proof (again, apart from the quantified version of Dirichlet's Theorem, due to de la Vall\'ee Poussin). If Schinzel's conditions are satisfied, the infinite product in~(\ref{eqn:C}) converges to a limit $C>0$. (See~\cite{AFG} for a proof, and for an interesting account of the background to the BHC.) Now the definite integral in~(\ref{eqn:E}) diverges to $+\infty$ as $x\to+\infty$, so $E(x)\to+\infty$ and hence, if the BHC is true, $Q(x)\to+\infty$ also, proving that there are infinitely many $t\in{\mathbb N}$ such that each $f_i(t)$ is prime.


\section{Heuristic argument for the BHC}\label{sec:heuristic}

According to the Prime Number Theorem, a good asymptotic estimate for the number $\pi(x)$ of primes $p\le x$ is obtained by integrating the probability of $t$ being prime, to give
\begin{equation}\label{eqn:better}
\pi(x)\sim\int_2^x\frac{dt}{\ln t}\quad\hbox{as}\quad x\to+\infty.
\end{equation}
(This is significantly more accurate than the widely-used estimate $x/\ln x$ for $\pi(x)$. Thus, for example, for the value of $\pi(x)$ for $x=10^{28}$, which may be 
found in \cite{OEIS}, entry A006880, the relative error of the estimate $x/\ln x$
is $-1.576\,\%$, while the relative error of the estimate (\ref{eqn:better}) is, approximately,
$10^{-12}\,\%$.)
This suggests that if $f(t)$ is a polynomial in ${\mathbb Z}[t]$ satisfying Bunyakovsky's conditions then the expected number of prime values of $f(t)$ for $t\le x$ might be estimated by
\begin{equation}\label{eqn-one-poly}
\int_a^x\frac{dt}{\ln f(t)}\quad\hbox{as}\quad x\to+\infty,
\end{equation}
where $a$ is chosen so that $f(t)>1$ for $t\ge a$. For a finite set of polynomials $f_1,\ldots, f_k$ satisfying Schinzel's conditions we can therefore make a first attempt at an estimate 
\begin{equation}\label{eqn-many-poly}
\int_a^x\frac{dt}{\prod_i\ln f_i(t)}
\end{equation}
for the expected number of $t\le x$ with all $f_i(t)$ prime, multiplying probabilities on the assumption that these polynomials behave independently of each other. However, they are not independent, and the Hardy--Littlewood constant $C$ is a product, over all primes $r$, of correction factors
\[\left(1-\frac{1}{r}\right)^{-k}\left(1-\frac{\omega_f(r)}{r}\right)\]
which replace the probability $\left(1-\frac{1}{r}\right)^k$ that $k$ randomly and independently chosen elements of ${\mathbb Z}_r$ are all non-zero with the probability that the product $f(t)=\prod_if_i(t)$ is non-zero mod~$(r)$, so that no $f_i(t)$ is divisible by $r$.

\begin{rema}
Note that the factor $C(f)$ modifies not only formula (\ref{eqn-many-poly}),
but also (\ref{eqn-one-poly}). As to (\ref{eqn:better}), we have $f(t)=t$, and it is easy
to see that in this case the Hardy--Littlewood constant is $C(f)=1$.
\end{rema}

Of course, the above is only a very brief outline of more persuasive heuristic 
arguments which can be used to support the BHC. See, however, the following
remark.

\begin{rema}
To the best of our knowledge, there is no completely adequate probabilistic model for the 
distribution of primes which would imply, even informally, the Bateman--Horn Conjecture (BHC).
Maybe, we must simply change direction and construct a model which would be based on
the BHC. At least, this latter conjecture has such strong supporting experimental data that a
model based on it just has no other way than to be an adequate description of a statistical
behaviour of primes. This task is not yet accomplished. 

However, we must keep in mind that the {\em density}, taken as a probability measure 
on $\mathbb{N}$ or on~$\mathbb{Z}$, is only finitely additive and not countable additive.
Therefore, it is impossible just to take a theorem of Probability Theory and apply it 
directly to the distribution of primes. The only particular case of the BHC which is  
proved up to now (the quantified version of Dirichlet's theorem due to 
de la Vall\'ee Poussin, mentioned above) was proved by analytic methods. As to Probability Theory, it is 
used as a heuristic tool, as a source of conjectures to be proved by other, non-probabilistic
methods.

In the absence of a proof, we have to rely on the fact that, time after time, the BHC 
produces estimates which agree remarkably closely with experimental data (see subsequent 
sections of this paper). Indeed, such evidence would regarded as convincing proof in 
other areas, such as Physics or Law.
\end{rema}


\section{Applying the BHC}\label{sec:applying}

We applied the BHC to the four triples of polynomials $f_i(t)$ in cases (a) to (d) of the six primes problem. In each case Maple can evaluate the definite integral in (\ref{eqn:E}) almost instantly, using numerical integration. (In the early1960s, when systems like Maple were unavailable and such programming had to be done `from scratch', Bateman and Horn simplified the integration by ignoring all coefficients of the polynomials $f_i$, replacing the definite integral in~(\ref{eqn:E}) with
\[\frac{1}{\prod_{i=1}^k\deg f_i}\int_2^x\frac{dt}{(\ln t)^k}.\]
The resulting estimates are asymptotically equivalent, but a little less accurate.

In case~(a) we have
\[f(t)=(12t+5)(3t+1)(2t+1).\]
This has three roots modulo each prime except $2$ and $3$, where it has one, so this gives us the terms in the infinite product~(\ref{eqn:C}). The product converges slowly, but by taking the product over all primes $r\le 10^9$ one can get a good approximation to the limit, namely $C(f)=5.71649719$. If we take $x=10^9$, for example,
 then the actual number of $t\le 10^9$ such that each $f_i(t)$ is prime is
$Q_a(10^9)=614\,423$, whereas the estimate is $E(x)=E_a(10^9)=615\,580.70$, representing 
a relative error of $+0.188\,\%$.
After the example $p=29$ given above, arising when $t=2$, the next primes in case~(a) appear when $t=14$, giving $p=173$, and then $t=26$, giving $p=317$.

The other three cases all give the same Hardy--Littlewood constant $C$ as in (a), leading to very similar BHC estimates. For instance, in case~(b) we have $Q_b(10^9)=615\,369$, while the estimate is $E_b(10^9)=615\,580.614$, an error of $+0.034\%$. In the other two cases $Q_c(10^9)=616\,509$ and $Q_d(10^9)=616\,289$, whereas the estimates are $E_c(10^9)=616\,720.62$ and $E_d(10^9)=616\,720.51$, giving errors of $+0.034\,\%$ and $+0.070\,\%$.

The accuracy of these estimates is comparable to that obtained in many other applications of the BHC, such as to twin or Sophie Germain primes (see also subsequent sections of this paper). Based on this evidence, we make the following conjecture: 
\begin{conj}
There are infinitely many primes $p$ satisfying the conditions in each of cases\/ $(a)$ to\/ $(d)$, and in particular there are infinitely many groups ${\rm PSL}_2(p)$ of order a product of six primes.
\end{conj}


\section{More than six primes}\label{sec:more}

One can easily modify the preceding arguments to obtain similar evidence for the existence of infinitely many groups ${\rm PSL}_2(p)$ in ${\mathcal S}_m$  for any $m\ge 6$. For instance, this can be done by replacing the conditions in case~(a) with
\[p-1=4r\quad\hbox{and}\quad p+1=2as\quad(a:=3^{m-5}),\]
where $p$, $r$ and $s$ are primes. Then $p\equiv 1$ mod~$(4)$ and $p\equiv -1$ mod~$(2a)$, so $p\equiv 2a-1$ mod~$(4a)$. Writing $p=4at+2a-1$ we have
\[r=\frac{p-1}{4}=at+\frac{a-1}{2}\quad\hbox{and}\quad s=\frac{p+1}{2a}=2t+1,\]
so we have an instance of the BHC with
\begin{equation}\label{eqn:three-poly}
f_i(t)=4at+2a-1,\quad at+\frac{a-1}{2}\quad\hbox{and}\quad 2t+1\quad(i=1,2,3).
\end{equation}

\begin{exam}
Let us take $m=7$, so that $a=3^{m-2}=9$. The polynomials we need to consider are
$$
f_1(t)=36t+17, \quad f_2(t)=9t+4 \quad \mbox{and} \quad f_3(t)=2t+1.
$$
For $r=2$ and $3$ the product $(36t+17)(9t+4)(2t+1)$ has a unique root modulo $r$, while for
all the other prime $r$ it has three roots. Therefore, the constant factor $C(f)$ remains the
same as in Section~\ref{sec:applying}. Then, taking, for example, $x=10^9$, we find that the
number of $t\le x$ for which the values of all the three polynomials are prime is
$Q(x)=556\,373$, while the BHC estimate gives $E(x)=556\,520.2$. The relative error of
this estimate is $0.026\,\%$.
\end{exam}

We make the following conjecture:
\begin{conj}
Given any $a=3^e$ ($e\ge 1$), the three polynomials 
in {\rm (\ref{eqn:three-poly})}
simultaneously take prime values for infinitely many $t\in{\mathbb N}$. In particular, 
given any $m\ge 6$ there are infinitely many groups ${\rm PSL}_2(p)$ of order a product 
of $m$ primes.
\end{conj}

This raises the question of whether one can {\em prove\/} that there is an infinite set of primes $p$ such that $\Omega(p^2-1)$ is bounded above, or more generally that there is an infinite set $\mathcal S$ of non-abelian finite simple groups $G$ with $\Omega(G)$ bounded above. Any such set $\mathcal S$ must contain only finitely many alternating groups, and hence any groups $G\in{\mathcal S}$ of Lie type must have bounded Lie rank, for otherwise their Weyl groups would involve alternating groups of unbounded degree.
Likewise, the field of definition ${\mathbb F}_q$ ($q=p^e$) of $G$ must have bounded degree $e$ over ${\mathbb F}_p$, for otherwise the order of the Sylow $p$-subgroups of $G$, which contain copies of the additive group of ${\mathbb F}_q$, would be an unbounded power of $p$. This suggests that the simplest way of constructing candidates for $\mathcal S$ is to use the groups ${\rm PSL}_2(p)$ as we have done here, with the Lie rank and the degree $e$ both equal to $1$.

\medskip

\begin{exam}\label{ex:unitary}
As another example where the rank and degree are both $1$, consider the simple unitary groups $G={\rm PSU}_3(p)$ ($=U_3(p)$ in ATLAS notation) of order $p^3(p^3+1)(p-1)/d$, $p>2$ prime,  where $d=\gcd(3,p+1)$. Since $p^3+1=(p+1)(p^2-p+1)$ it is not hard to see that $\Omega(G)\ge 9$, attained if and only if $\Omega(p^2-1)=6$ (as in the case of ${\rm PSL}_2(p)$) and  $p^2-p+1$ or $(p^2-p+1)/3$ is prime. (The latter case arises when $p\equiv -1$ mod~$(3)$, so that $d=3$.) When we considered ${\rm PSL}_2(p)$ in Sections~\ref{sec:factors} and \ref{sec:applying} the prime $p$ was represented by $f_1(t)$ in each of cases (a) to (d), so
one can apply the BHC to ${\rm PSU}_3(p)$ by adding the polynomial $f_4(t)=f_1(t)^2-f_1(t)+1$ to the triples $f_i(t)$ ($i=1,2,3$) used earlier in cases (b) and (d). In cases~(a) and~(c) we have $f_1(t)\equiv -1$ mod~$(3)$ so $f_1(t)^2-f_1(t)+1$  is divisible by $3$ and we can take 
$f_4(t)=(f_1(t)^2-f_1(t)+1)/3$. The following example gives strong evidence that there are infinitely many groups ${\rm PSU}_3(p)\in{\mathcal S}_9$.
\end{exam}

\begin{exam}\label{ex:four-poly}
Let us consider one particular example coming from case~(c): we have
$$
f_1(t)=12t-1, \quad f_2(t)=6t-1, \quad f_3(t)=t, \quad \mbox{and} \quad
f_4(t)=\frac{f_1(t)^2-f_1(t)+1}{3}=48t^2-12t+1.
$$
In order to compute the constant $C(f)$ we need to know the value $\omega_f(r)$
which is the number of roots of the product $f(t)=f_1f_2f_3f_4$ in $\mathbb{Z}_r$, with
$r$ prime. For $r=2$ and $3$ this product has a single root. For $r>3$ the factors
$f_1,f_2,f_3$ provide us with three roots, while the quadratic polynomial $f_4$ may
have two roots or no roots. Let us look at this case in more detail.

The discriminant of $f_4$ is $-48=4^2\cdot(-3)$. Therefore, $f_4$ has two roots 
in $\mathbb{Z}_r$ if $-3$ is a quadratic residue modulo $(r)$, and no roots if not.
Recall the notation for the Legendre symbol: for $p$ prime and $q\in\mathbb{Z}$
$$
\left(\frac{q}{p}\right)=
\begin{cases}
0\quad\; \hbox{if $q\equiv 0$ mod~$(p)$;}\\
1\quad\; \hbox{if $q$ is a quadratic residue mod~$(p)$;}\\
-1 \;\; \hbox{otherwise.}
\end{cases}
$$
(See~\cite[Chapter~7]{JJ} for quadratic residues and the Legendre symbol.) This symbol is
multiplicative: for all $q, q'\in{\mathbb Z}$ we have
$$
\left(\frac{qq'}{p}\right)=\left(\frac{q}{p}\right)\left(\frac{q'}{p}\right),
\quad \mbox{therefore} \quad 
\left(\frac{-3}{p}\right)=\left(\frac{-1}{p}\right)\left(\frac{3}{p}\right).
$$
The value of $\left(\frac{-1}{p}\right)$ is known: it is 1 if $p=4k+1$,
and $-1$ if $p=4k-1$. As to $\left(\frac{3}{p}\right)$, we may use
Gauss' Law of Quadratic Reciprocity: for $p,q$ prime we have 
$\left(\frac{q}{p}\right)=\left(\frac{p}{q}\right)$ if one or both $p$ and $q$ are
of the form $4k+1$, and $\left(\frac{q}{p}\right)=-\left(\frac{p}{q}\right)$ if both
$p$ and $q$ are of the form $4k-1$. Collecting all these data, we find out that
$$
\left(\frac{-3}{p}\right) = \left(\frac{p}{3}\right) =
\left\lbrace 
\begin{array}{cl} 
1  & \mbox{if } p \equiv \; 1 \; \mbox{ mod } 3, \\
-1 & \mbox{if } p \equiv -1 \mbox{ mod } 3.
\end{array}
\right.
$$

We may now conclude: for a prime $r>3$ the product $f(t)=f_1(t)f_2(t)f_3(t)f_4(t)$
has five roots in~$\mathbb{Z}_r$ if $r \equiv  1 \mbox{ mod } (3)$, and has three roots
if $r \equiv  -1 \mbox{ mod } (3)$. This permits us to compute the constant $C(f)$:
taking the product over the primes $r\le 10^9$ we get $C(f)=12.10128533$. The number of
$t\le 10^9$ such that all $f_i(t)$ are prime is $Q(10^9)=30\,452$; the estimate of
this number given by the BHC is $30\,504.71$; the relative error is $0.173\,\%$.
\end{exam}


\section{Permutation groups of prime degree}\label{sec:projective}

Although we have concentrated in this paper on the six primes problem, our involvement with the BHC began with a different problem in group theory. This was motivated by our study of Klein's papers~\cite{Klein7,Klein11} on equations (and hence coverings and permutation groups) of degree $7$ and $11$, and our desire to extend his results to all primes (see~\cite{JZ1} for details). In doing this we encountered an important but rarely discussed gap in the classification of transitive permutation groups of prime degree, a problem going back two and a half centuries to Lagrange. Building on early work by Galois and Burnside, the classification 
of finite simple groups implies that these groups are as follows: 

\begin{itemize}
\item[(a)]	subgroups of ${\rm AGL}_1(p)$ containing the translation group, for primes $p$;
\item[(b)]	alternating and symmetric groups ${\rm A}_p$ and ${\rm S}_p$, for primes $p\ge 5$;
\item[(c)]	${\rm PSL}_2(11)$ and Mathieu groups  ${\rm M}_{11}$ and ${\rm M}_{23}$, of 
			degrees $11, 11$ and $23$;
\item[(d)]	subgroups $G$ of ${\rm P\Gamma L}_n(q)$ containing ${\rm PSL}_n(q)$, in those 
			cases when the natural degree \linebreak $d=(q^n-1)/(q-1)$ is prime.
\end{itemize}

In (c), ${\rm PSL}_2(11)$ has two actions of degree $11$, on the cosets of two conjugacy classes of subgroups isomorphic to ${\rm A}_5$; the Mathieu groups act on the points of the Steiner systems $S(4,5,11)$ and $S(4,7,23)$. In (d), $G$ acts on the $d$ points (and also the $d$ hyperplanes if $n\ge 3$) of the projective geometry ${\mathbb P}^{n-1}({\mathbb F}_q)$ of dimension $n-1$ over 
${\mathbb F}_q$.

It is an open problem whether the degree $d$ in case~(d) is prime for infinitely many pairs $(n,q)$.
Such {\sl projective primes}, as we call them, include the Fermat primes, of the form $q+1=2^{2^f}+1$, for $n=2$, and the Mersenne primes, of the form $2^n-1$, for $q=2$. Five Fermat primes are known, and it is conjectured that there are no more; at the time of writing $51$ Mersenne primes are known, and it is conjectured that there are infinitely many of them. In investigating this problem (see~\cite{JZ2} for details), having nothing new to contribute to the extensive research on those very difficult topics, we restricted our attention to the cases $n, q\ge 3$.

If we write $q=p^e$, we are asking whether $p$ and
\begin{equation}
d:=p^{(n-1)e}+p^{(n-2)e}+\cdots+p^e+1
\end{equation}
can both be prime for infinitely many $p$. Clearly, this requires $n$ to be prime, and a simple argument involving cyclotomic polynomials shows that $e$ must be a power of $n$, possibly equal to $n^0=1$.

We concentrated on the simplest and most frequently-arising case $n=3$, $e=1$, applying the BHC to the polynomials
\[f_1(t)=t,\quad f_2(t)=t^2+t+1.\]
Typical results obtained were
\begin{itemize}
\item	$E(10^{10})=1.579642126 \times 10^7$ and $Q(10^{10})=15\, 801\, 827$, an error of 
		$-0.03420956\,\%$, 
\item 	$E(10^{11})=1.292974079 \times 10^8$ and $Q(10^{11})=129\, 294\, 308$, an error of 
		$+0.00239757\,\%$.
\end{itemize}
For other pairs $(n,e)$,  such as $(5,1)$ and $(3,3)$, the results were good but much less convincing, since the primes increase so rapidly that relatively few of them were within our computing range. Nevertheless, on the basis of this evidence we make the following conjecture:

\begin{conj}
For each prime $n\ge 3$ there are infinitely many prime powers $q$ such that ${\rm PSL}_n(q)$ is a transitive permutation group of prime degree $d=(q^n-1)/(q-1)$.
\end{conj}


\section{Linear groups of prime degree}\label{sec:linear}

In~\cite{DZ} Dixon and Zalesskii classified the irreducible finite subgroups $G\le{\rm SL}_d({\mathbb C})$ of prime degree ~$d$, in the case where the socle $S$ of their image $\overline G$ in ${\rm PSL}_d({\mathbb C})$ is non-abelian and acts primitively, that is, preserving no non-trivial direct sum decomposition of ${\mathbb C}^d$. (In this section  `degree' and `primitive' always refer to the dimension and structure of the vector space on which $G$ acts, and not to any permutation representation of $G$.) In this case $S$ is simple and $\overline G\le{\rm Aut}\,S$. Theorem~1.2 of~\cite{DZ} gives a finite list of families of simple groups $S$ which can arise, with necessary and sufficient conditions on $d$ and their parameters for such groups $G$ to exist. It is unknown whether some of these families are finite or infinite. Here we give a brief outline of how we have used the BHC to give strong evidence that they are infinite (see~\cite{JZ2} for full details).

A typical example has $S\cong {\rm PSU}_n(q)$ for $n\ge 3$, where $q$ is a prime power $p^e$ and the degree
\[d=\frac{q^n+1}{q+1}=q^{n-1}-q^{n-2}+\cdots-q+1\]
is prime. The pair of polynomials
\begin{equation}\label{eq:PSU}
f_1(t)=t\quad\hbox{and}\quad f_2(t)=t^{(n-1)e}-t^{(n-2)e}+\cdots-t^e+1
\end{equation}
satisfy the conditions of Schinzel's Hypothesis if $n$ is prime and $e$ is a power of $n$.
This is the same as the condition in Section~\ref{sec:projective} for $(q^n-1)/(q-1)$ to be a projective prime, and indeed the substitution $t\mapsto -t$ shows that the Hardy--Littlewood constants $C(f)$ are the same in both cases: they are equal to 1.521730.

We may compare the results, for example, for $x=10^{10}$. The definite integrals 
in~(\ref{eqn:E}) are very similar: in fact, they differ by 0.456 (and the estimates $E(x)$ 
thus differ by $C(f)\cdot 0.456\approx 0.7$). Therefore, we may suppose that the actual numbers
$Q(x)$ will also be close to each other. And, indeed, we have seen above that the value
of $Q(x)$ for $f_1(t)=t$, $f_2(t)=t^2+t+1$ was $15\,801\,827$, while for the polynomials
$f_1(t)=t$, $f_2(t)=t^2-t+1$ we have $15\,801\,414$ (the relative error in this case is
$-0.0316\,\%$). It is interesting to note that, again for $x=10^{10}$, the computation of
the estimate $E(x)$ on a particular (modest) laptop we have used took 
0.002~seconds, while the computation of the exact value of $Q(x)$ took 54 hours. 
The computation of $E(x)$ for $x=10^{30}$ {\em with\/ $30$~digits of accuracy}\/
takes 2.4 seconds, while the computation of $Q(x)$ for this value of $x$ is far 
beyond our reach. Regrettably, we don't have a (conjectural) upper bound for the error term.

Another family appearing in~\cite{DZ} consists of groups $S\cong{\rm PSL}_2(q)$ where $q$ and the degree $d=(q-1)/2$ are both prime. This is equivalent to $d$ being a Sophie Germain prime, a case where the BHC has already provided strong evidence of infinitely many examples.

A more difficult family in~\cite{DZ} has $S\cong{\rm PSL}_2(q)$ for prime degrees 
$d=\frac{q+1}{2}$, where $q=p^{2^k}\ge 5$ for an odd prime $p$ and integer $k\ge 0$. 
To apply the BHC we took
\[f_1(t)=2t+1\;(=p)\quad\hbox{and}\quad f_2(t)=\frac{(2t+1)^{2^k}+1}{2}=\sum_{i=1}^{2^k}\binom{2^k}{i}2^{i-1}t^i+1\;(=d)\]
for some fixed $k\ge 0$. In this case the calculation of the Hardy--Littlewood constants $C(f)$ for $k\ge 1$ is less straightforward (the case $k=0$ is similar to the Sophie Germain primes problem). We have
\[
\omega_f(r) \,=\, \left\lbrace
\begin{array}{ll}
0   & \mbox{\rm if } r=2, \\
2^k+1 & \mbox{\rm if } r\equiv 1\; \mbox{\rm mod } (2^{k+1}), \\
1	& \mbox{\rm otherwise},
\end{array}
\right.
\]
leading to the results for small $k$ shown in Table~\ref{tab:d=2^k}:

\begin{table}[htbp]
\begin{center}
\begin{tabular}{c|c|c|c|c}
$k$  	& $C(f)$      & $Q(10^9)$   & $E(10^9)$      & relative error  \\
\hline
1       & $4.426783$  & 5\,448\,994 & 5\,448\,648.05 & $-0.006$\,\%    \\
2		& $10.433814$ & 6\,373\,197 & 6\,365\,668.39 & $-0.118$\,\%    \\
3		& $7.885346$  & 2\,394\,012 & 2\,395\,075.38 & $0.044$\,\%		\\
4		& $14.642571$ & 2\,219\,445 & 2\,218\,975.66 & $-0.021$\,\%
\end{tabular}
\end{center}
\smallskip
\caption{$Q(10^9)$ is the number of $t\le 10^9$ such that both $f_1(t)$ and $f_2(t)$
are prime; $E(10^9)$ is the BHC estimate for $Q(10^9)$.}
\label{tab:d=2^k}
\end{table}

On the basis of this we conjecture that for each $k\ge 0$ there are infinitely many primes $p$ such that $d$ is prime.

The above example is, in fact, the particular case $n=1$ of a more general family involving the symplectic groups $S={\rm PSp}_{2n}(q)$ where $d=(q^n+1)/2$ is prime, $n=2^j$ for some integer $j\ge 0$, and $q$ is as before, but without the restriction $q\ge 5 $ when $n>1$. 
(Note that ${\rm PSp}_2(q)\cong{\rm PSL}_2(q)$.) For a fixed pair $j, k$ we applied the BHC to the polynomials
\[f_1(t)=2t+1\;(=p)\quad\hbox{and}\quad f_2(t)=\frac{(2t+1)^{2^{j+k}}+1}{2}\;(=d).\]
These are the same as the preceding pair $f_1, f_2$, but with $j+k$ replacing $k$, so the same estimates $E(x)$ and search results $Q(x)$ apply in this case.

There are several other potentially infinite   families of groups in~\cite[Theorem~1.2]{DZ}, but since they involve exponential functions rather than polynomials they are beyond the scope of the BHC.

In a corrigendum to~\cite{DZ} Dixon and Zalesskii showed that if $G$ is primitive, and the socle $S$ is non-abelian and imprimitive, then $G$ has an imprimitive commutator subgroup $G'\cong{\rm PSL}_n(q)$ where $d=(q^n-1)/(q-1)$, with $q$ odd or $q=2$. Our results on projective primes apply in this case for odd $q$, and they are also relevant to~\cite{DZ04} where the same authors have considered imprimitive linear groups of prime degree, transitively permuting the one-dimensional subspaces in a direct sum decomposition.


\section{Other group-theoretic applications of the BHC}\label{sec:other}

Cameron, Manna and Mehatari~\cite{CMM} have recently studied the {\em power graph\/} $P(G)$ of a finite group~$G$. The vertices of $P(G)$ are the elements of $G$, and a pair of them are joined by an edge if one of them is a power of the other. In their Theorem~1.3 they characterise those non-abelian finite simple groups $G$ for which $P(G)$ is a {\em cograph}. (This is a graph with no induced subgraph isomorphic to the path $P_4$ with four vertices; cographs form the closure of the one-vertex graph $K_1$ under the operations of disjoint union and complement.) In their Problem~1.4 they ask whether there are infinitely many such groups $G$.
 
One family appearing in their characterisation consists of the groups ${\rm PSL}_2(q)$ for odd prime powers $q\ge 5$ such that $(q\pm 1)/2$ are each either a prime power or a product of two primes. The groups ${\rm PSL}_2(p)$ in cases~(a) and~(b) of the six primes problem (see Section~\ref{sec:factors}) satisfy this condition, so the BHC gives strong evidence of a positive answer to their question.
 
In~\cite{ADP} Amarra, Devillers and Praeger have recently constructed families of block designs which have interesting symmetry properties (a group of automorphisms acting transitively on blocks, and transitively but imprimitively on points) and which maximise various parameters. Their constructions, based on finite fields, require certain quadratic polynomials to take prime power values. By using the BHC, together with some extensions from primes to prime powers, we have in~\cite{JZ3} provided strong evidence that these families are infinite. 


\section{Acknowledgements}

The authors are grateful to Peter Cameron and Cheryl Praeger for helpful comments on applying the BHC to their work. Alexander Zvonkin is supported by the French ANR project {\sc Combin\'e} (ANR-19-CE48-0011).
  

\end{document}